\documentclass[twocolumn]{autart}    % Enable this line and disable the
\usepackage{amssymb,amsmath}
\usepackage{graphicx}
\usepackage{color}
\usepackage{subfigure}        % Include this line if your
\emergencystretch=\maxdimen
\newcommand{\col}{$\upshape{col}$}

\newtheorem{rmk}{Remark}

\begin{document}

\begin{frontmatter}

\title{Distributed event-triggered aggregative optimization with applications to  price-based energy management\thanksref{footnoteinfo}} % Title, preferably not more

\thanks[footnoteinfo]{This work was supported  in part by the National Natural Science Foundation of China (NSFC, Grant Nos. 62273145 and 62303394), in part by the Beijing
Natural Science Foundation (Grant No. 4222053), and in part by the Natural Science Foundation of Xinjiang Uygur Autonomous Region (Grant No. 2022D01C694). (Corresponding author: Feng~Xiao.) }

\author[Beijing1,Beijing2,Urumqi]{Xin Cai}\ead{caixin\_xd@126.com},    % Add the
\author[Beijing1,Beijing2]{Feng Xiao}\ead{fengxiao@ncepu.edu.cn},               % e-mail address
\author[Beijing1]{Bo Wei}\ead{bowei@ncepu.edu.cn},
\author[Beijing3]{Aiping Wang}\ead{aiping@ncepu.edu.cn}
 % Add the
%\author[Urumqi]{Xinyuan Nan}\ead{xynan@xju.edu.cn},               % e-mail address
%\author[Urumqi]{Binpeng Gao}\ead{gbp\_xd@sina.com},  % (ead) as shown
%\author[Beijing]{Mei Yu}\ead{meiyu@ncepu.edu.cn},
%\author[Beijing]{Fang Fang}\ead{ffang@ncepu.edu.cn}
%
\address[Beijing1]{ State Key Laboratory of Alternate Electrical
Power System with Renewable Energy Sources, North China Electric Power University, Beijing, 102206, China}

\address[Beijing2]{School of Control and Computer Engineering, North China Electric Power University, Beijing, 102206, China}  %
\address[Urumqi]{School of Electrical Engineering, Xinjiang University, Urumqi, 830017, China}             % full addresses
    % here.
\address[Beijing3]{School of Mathematics and Physics, North China Electric Power University, Beijing, 102206, China}

                             % keyword list or with the
                                          % help of the Automatica
                                          % keyword wizard

\begin{abstract}
This paper studies a distributed continuous-time aggregative optimization problem, which is a fundamental problem in the price-based energy management. The objective of the distributed aggregative optimization is to minimize the sum of local objective functions, which have a specific expression that relies on agents' own decisions and the aggregation of all agents' decisions. To solve the problem, a novel distributed continuous-time algorithm is proposed by combining gradient dynamics with a dynamic average consensus estimator in a two-time scale. The exponential convergence of the proposed algorithm is established under the assumption of a convex global cost function by virtue of the stability theory of singular perturbation systems. Motivated by practical applications, the implementation of the continuous-time algorithm with event-triggered communication is investigated.
Simulations on the price-based energy management of distributed energy resources are given to illustrate the proposed method.
\end{abstract}

\end{frontmatter}

\section{Introduction}
In the last decade, steering a networked dynamical system to an optimal steady state has been a hot topic. In some engineering scenarios, the optimal steady states have appropriate economic interpretations. For example, in the economic dispatch of power systems, price-anticipating distributed energy resources (DERs) dynamically adjust their generations in a cooperative manner to minimize the sum of local optimization objectives. The local optimization functions depend on both local cost of each DER and the pricing function determined by total power generations \cite{Cai.2021,Cai.2022,Gharesifard.2016,Ye.2017,Ye.2021}, and are described by $f_i(P_i,\sigma(P))$$=$$c_i(P_i)-p(\sigma(P))P_i$, $\forall i\in\{1,\ldots,N\}$,
where $P_i$ denotes the energy generation of the $i$-th DER, $c_i(P_i)=a_iP_i^2+b_iP_i+d_i$ is the generation cost with parameters $a_i$, $b_i$, and $d_i$, the electricity price $p(\sigma)=200-0.1N\sigma(P)$ with $\sigma(P)=(1/N)\sum_{j=1}^NP_j$ and $P=[P_1,\ldots,P_N]^T$. Thus, the global cost minimization problem, solved by the network of DERs, is given by
\begin{equation} \label{ex}
\min_{P\in\mathbb{R}^N} \sum_{i=1}^N f_i(P_i,\sigma(P)).
\end{equation}
In addition, optimization functions in similar forms have appeared in the network congestion control \cite{Grammatico.2017} and coverage control of sensor networks \cite{Romano.2020}. However, individual costs were minimized by distributed algorithms proposed in \cite{Cai.2021,Cai.2022,Gharesifard.2016,Grammatico.2017,Romano.2020,Ye.2017,Ye.2021} to obtain local optima. In the distributed aggregative optimization \cite{Li.2021}, like \eqref{ex}, the global optimization goal is a sum of local objective functions relying on agents' own decision variables and  an aggregation of all agents' decisions.

As a special class of distributed optimization problems, the aggregative optimization problem was first proposed in \cite{Li.2021}. In contrast to the distributed optimization problem in which each agent estimates global optimal solutions by interchanging information containing its own decision, the aim of the distributed aggregative optimization is to optimize each agent's own decision in the global sense. Moreover, the interchanged information only includes estimations of the aggregation, which protects privacy of agents' decisions. Recently, the distributed aggregative optimization with local constraints was studied in \cite{Wang.2022}. Taken optimization functions and local constraints changing over time into consideration, a distributed online aggregative optimization method was proposed in \cite{Carnevale.2022a}. Different from the discrete-time system studied in \cite{Carnevale.2022a,Li.2021,Wang.2022}, a network of continuous-time dynamical systems is considered in this paper. To find the global optimum of an aggregative optimization problem, a distributed algorithm based on local information is proposed. From the viewpoint of control theory, the designed distributed algorithm can be regarded as a feedback control system,  whose stability is analyzed by the Lyapunov stability theory. The main contributions of this paper are summarized as follows.

(i) A novel distributed continuous-time algorithm is designed to solve an aggregative optimization problem. Combined with a dynamic average consensus estimator for the aggregation information, a distributed gradient based algorithm is proposed in a two-time scale. The estimator is updated in a fast time scale and the decision of each agent evolves in a slow one. Different from linear convergence rate obtained in \cite{Carnevale.2022a,Li.2021,Wang.2022}, the exponential convergence of the designed algorithm is guaranteed under the assumption of a convex global cost function and is proved by the Lyapunov stability theory and the singular perturbation theory.

(ii) To reduce communication loads caused by continuous-time communication, a distributed and asynchronous event-triggered scheme is proposed based on each agent's own information rather than neighbors' last broadcast information, which is required in the event-triggered communication schemes designed in \cite{Cai.2021a,Chen.2016,Dai.2020,Liu.2020,Wang.2019a,Yu.2021}. The exponential convergence is preserved and the Zeno behavior, which is a phenomenon that an infinite number of events occur in a bounded time interval, is excluded.

(iii) In contrast to non-cooperative DERs scenarios \cite{Gharesifard.2016,Ye.2017}, a cooperative scenario, in which DERs cooperate to minimize the sum of local cost functions coupling with a price function, is presented as an example. The proposed scenario is characterized by the distributed aggregative optimization, thus it cannot be addressed by using the  methods in \cite{Gharesifard.2016,Ye.2017}.

The rest of this paper is organized as follows. In Section \ref{sec2}, the problem formulation is given. In Section \ref{sec3}, a distributed continuous-time algorithm based on an event-triggered broadcasting scheme is designed.  A numerical example is provided in Section \ref{sec4}. Finally, some conclusions are stated in Section \ref{sec5}.

Notations: $\mathbb{R}$ and $\mathbb{R}_{\geq0}$ denote the sets of real and non-negative real numbers, respectively. $\mathbb{R}^n$ is the $n$-dimensional real vector space. $\mathbb{R}^{n\times m}$ denotes the set of $n \times m$ real matrices. Given a vector $x\in\mathbb{R}^n$, $\|x\|$ is its Euclidean norm. $A^T$ and $\|A\|$ are the transpose and the spectral norm of matrix $A\in \mathbb{R}^{n\times n}$, respectively. For matrices $A$ and $B$, $A\otimes B$ denotes their Kronecker product. Let $\col(x_1,\ldots,x_n)=[x_1^T,\ldots,x_n^T]^T$.  $1_n$ and $0_n$ are $n$-dimensional column vectors with all elements being ones and zeros, respectively. $I_n$ denotes an $n\times n$ identity matrix. $B_r$ denotes the ball $\{x\in\mathbb{R}^n| \|x\|\leq r\}$ centered at the origin. $\nabla f(x)$ is the gradient of function $f(x)$ with respect to $x$.

\section{Problem formulation} \label{sec2}
Consider a network of agents which are indexed by the set $\mathcal{I}=\{1,\ldots,N\}$ and communicate with each other via an undirected and connected graph $\mathcal{G}$, which is not weighted. Here, it is assumed that agents are with continuous-time single-integrator dynamics
\begin{equation} \label{dyn}
\dot{x}_i=u_i,\ \ i\in\mathcal{I}
\end{equation}
where $x_i\in\mathbb{R}^{n_i}$ is the decision variable and $u_i\in\mathbb{R}^{n_i}$ is the control input of agent $i$. Motivated by the formation control of networked agents \cite{Li.2021} and energy management of DERs \cite{Gharesifard.2016}, assign agent $i$ a local cost function $f_i(x_i,\sigma(x)):\mathbb{R}^{n}\rightarrow \mathbb{R}$ with $n=\sum_{i=1}^N n_i$, where $\sigma(x): \mathbb{R}^n\rightarrow \mathbb{R}^m$ is an aggregative function relying on the decisions of all agents and is defined by $\sigma(x)=(1/N)\sum_{i=1}^N\phi_i(x_i)$. Function $\phi_i:\mathbb{R}^{n_i}\rightarrow \mathbb{R}^m$ is only known to agent $i$. For agent $i$ with dynamics \eqref{dyn}, the objective is to design $u_i$ by using only local information so that all agents' decisions $x$ arrive at the optimal decision $x^*=\col(x_1^*,\ldots,x_N^*)$, which is the minimum of the following aggregative optimization problem
\begin{equation} \label{op}
\min_{x\in\mathbb{R}^n} \sum_{i=1}^Nf_i(x_i,\sigma(x)).
\end{equation}

Define the global cost function $f(x)$$=$$\sum_{i=1}^Nf_i(x_i,\sigma(x))$. The following assumptions ensure the solvability of problem \eqref{op}.

\begin{assum} \label{as2}
The convex function $f(x)$ is continuously differentiable.
\end{assum}

\begin{assum} \label{as3}
The map $\nabla f(x)$ is $l$-locally Lipschitz continuous with $l>0$, and there exist a positive constant $\kappa$ and an open set $S\supset x^*$ such that $\|\nabla f(x)\|^2\geq (1/\kappa)\|x-x^*\|^2$, $\forall x\in S$.
\end{assum}

Assumption \ref{as3} is similar to the metrically subregular given in \cite{Liang.2019}, which relaxes the requirement of strongly convex functions to ensure the exponential convergence.
Under Assumption \ref{as2}, problem \eqref{op} is a convex optimization problem. Based on the first-order optimal condition in \cite{Boyd.2004}, a necessary and sufficient condition of optimal decision $x^*$ is given by  $\partial f_i(x_i^*,\sigma(x^*))/\partial x_i+(\nabla \phi_i(x_i)/N)(\sum_{j=1}^N\partial f_j(x_j^*,\sigma(x^*))/\partial \sigma)$$=$$0_{n_i}$ for all $i\in\mathcal{I}$. According the optimal condition, a centralized optimization algorithm is given by
\begin{equation} \label{ci}
\dot{x}_i\!=\!-\frac{\partial f_i(x_i,\sigma(x))}{\partial x_i}\!-\!\frac{\nabla \phi_i(x_i)}{N}\!\sum_{j=1}^N\!\frac{\partial f_j(x_j,\sigma(x))}{\partial \sigma}.
\end{equation}

\begin{lem} \label{thm1}
Suppose that Assumptions \ref{as2}-\ref{as3} hold. Algorithm \eqref{ci} converges to the optimal decision $x^*$ at an exponential convergence rate no less than $\kappa/(1+2l)$.
\end{lem}

Lemma \ref{thm1} can be easily obtained by the analytical methods given in \cite{Liang.2019,Yi.2020}. The proof of Lemma \ref{thm1} is omitted for the limited space. Note that it is seen from \eqref{ci} that each agent needs to know the information of aggregative function $\sigma(x)$ and $\sum_{j=1}^N\partial f_j/\partial \sigma$ by a coordinator which knows the global objective, receives all agents' decisions and sends information to each agent. To make algorithm \eqref{ci} more reliable and robust, in the subsequent section, algorithm \eqref{ci} will be exploited to propose a distributed algorithm for solving aggregative optimization problem \eqref{op}.

\begin{rmk}
In aggregative games \cite{Carnevale.2022,Koshal.2016,Liang.2017,Liang.2022,Shakarami.2022}, the individual optimization objectives have a similar form of $f_i(x_i,\sigma(x))$. Essentially, different from the local optima obtained by distributed Nash equilibrium seeking algorithms designed in  \cite{Carnevale.2022,Koshal.2016,Liang.2017,Liang.2022,Shakarami.2022}, the global optimum of problem \eqref{op} is achieved in this paper.
\end{rmk}

\section{Main Results} \label{sec3}
In this section, a distributed continuous-time algorithm with an event-triggered broadcasting scheme is designed for agents with dynamics \eqref{dyn} to achieve their own optimal decisions of aggregative optimization problem \eqref{op}.

In a distributed setting, each agent is required to estimate some global information by a dynamic average consensus estimator. For agent $i$, define $\eta_{i1}\in\mathbb{R}^m$ and $\eta_{i2}\in\mathbb{R}^m$ as the estimations of aggregative function $\sigma(x)$ and  $(1/N)\sum_{j=1}^N\partial f_j/\partial \sigma$, respectively. Denote $\eta_i$$=$$\col(\eta_{i1},\eta_{i2})$ and introduce $w_i\in\mathbb{R}^{2m}$ as an auxiliary variable. Let $\Theta_i$$=$$\col(\phi_i(x_i),\partial f_i/\partial \sigma)$. The distributed algorithm is designed by
\begin{equation} \label{alg}
\dot{x}_i=-\frac{\partial f_i(x_i,\eta_{i1})}{\partial x_i}-\eta_{i2}\nabla \phi_i(x_i)
\end{equation}
with a dynamic average consensus estimator, given by
\begin{equation} \label{dacd}
\begin{aligned}
\delta \dot{\eta}_i&=-\eta_i-\sum_{j\in\mathcal{N}_i} (\hat{\eta}_i-\hat{\eta}_j)-\sum_{j\in\mathcal{N}_i} (\hat{w}_i-\hat{w}_j)\!+\!\Theta_i,\\
\delta \dot{w}_i&=\sum_{j\in\mathcal{N}_i}(\hat{\eta}_i-\hat{\eta}_j).
\end{aligned}
\end{equation}
In \eqref{dacd}, $\mathcal{N}_i$ is the set of agent $i$'s neighbors and $\delta$ is a positive perturbation parameter. $\hat{\eta}_i(t)=\eta_i(t_{k_i})$ and $\hat{w}_i(t)=w_i(t_{k_i})$ for $t\in[t_{k_i},t_{k_i+1})$ are the broadcast information to the neighboring agents of agent $i$ at time $t_{k_i}, k_i=1,2,...$. Let $e_i=\col(\hat{\eta}_i-\eta_i,\hat{w}_i-w_i)$. The triggering instant is determined by the following rule:
\begin{equation} \label{tc}
t_{k_i+1}=\inf \{t>t_{k_i}| \  \|e_i\|\geq\beta_{i1}e^{-\beta_{i2} t}\},
\end{equation}
where  $\beta_{i1}, \beta_{i2}>0$. The overall control framework for the multi-agent system is shown in Fig.~1.

\begin{figure}[!t]
  \centering
  \includegraphics[width=6cm]{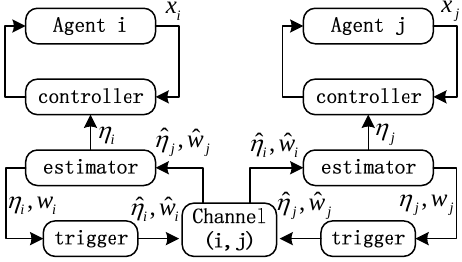}
  \caption{The control framework for the multi-agent system.}
  \label{fig8}
\end{figure}

\begin{rmk}
The event-triggered frequency in \eqref{tc} greatly depends on parameters $\beta_{i1}$ and $\beta_{i2}$ \cite{Seyboth.2013}. To ensure positive inter-event intervals, $\beta_{i2}$ should be selected suitably such that state $(\eta_i,w_i)$ converges faster than $\beta_{i1}e^{-\beta_{i2} t}$ decays. In addition, for fixed $\beta_{i2}$, if $\beta_{i1}$ is chosen to be large, the inter-event intervals become large.
\end{rmk}

Laplacian matrix of graph $\mathcal{G}$ is denoted by $L$. Denote $\eta$ $=$ $\col(\eta_1,\ldots,\eta_N)$, $w$ $=$ $\col(w_1,\ldots,w_N)$, $F(x,\eta_1)$ $=$ $\col(\partial f_1(x_1,\eta_{11})/\partial x_1,\ldots,\partial f_N(x_N,\eta_{N1})/\partial x_N)$, $\eta_2$ $=$ $\col(\eta_{12},\ldots,\eta_{N2})$, $\nabla \phi(x)$$=$$\col(\nabla \phi_1(x_1),\ldots, \nabla \phi_N(x_N))$, $\Theta$$=$$\col(\Theta_1,\ldots,\Theta_N)$, $E_1=\col(\hat{\eta}_1-\eta_1,\ldots,\hat{\eta}_N-\eta_N)$ and $E_2=\col(\hat{w}_1-w_1,\ldots,\hat{w}_N-w_N)$. The compact form of system \eqref{alg}-\eqref{dacd} is given by
\begin{equation} \label{cf}
\begin{aligned}
\dot{x}&=-F(x,\eta_1)-\eta_2\nabla \phi(x),\\
\delta \dot{\eta}&=-\eta-(L\otimes I_{2m})\eta-(L\otimes I_{2m})w+\Theta\\
&\ \ \ -(L\otimes I_{2m})(E_1+E_2),\\
\delta \dot{w}&=(L\otimes I_{2m})\eta+(L\otimes I_{2m})E_1.
\end{aligned}
\end{equation}

\begin{rmk}
A significant difference between the algorithms proposed in \cite{Carnevale.2022a,Li.2021,Wang.2022} and our designed algorithm is in the problem setup. This paper studies the distributed aggregative optimization problem among continuous-time agents with event-triggered communication, unlike the work \cite{Carnevale.2022a,Li.2021,Wang.2022}, in which the problem among discrete-time agents with periodic communication was studied. Furthermore, different from proportional dynamic average consensus estimators for the estimation of certain global information \cite{Carnevale.2022a,Li.2021,Wang.2022}, a proportional-integral dynamic average consensus estimator (7) is used in this paper to ensure accurate estimation. In addition, this paper relaxes the requirement of strongly convex objective functions given in \cite{Carnevale.2022a,Li.2021}, just convex objective functions are required to ensure the exponential convergence of algorithm (6)-(7).
\end{rmk}

\begin{lem} \label{lem2}
$x^*$ is an optimal decision of problem \eqref{op} if and only if $(x^*,\eta^*,w^*)$ is an equilibrium of system \eqref{cf} with $\eta_1^*=\cdots=\eta_N^*=(1/N)\sum_{j=1}^N\Theta_j$.
\end{lem}

\textbf{Proof.} \textit{Necessity}: Assume that $(x^*,\eta^*,w^*)$ is an equilibrium of system \eqref{cf}. At the equilibrium point, the measurement errors $E_1$ and $E_2$ are  zeros. We have that
\begin{subequations}
\begin{align}
0_n&=-F(x^*,\eta_1^*)-\eta_2^*\nabla \phi(x^*), \label{cf11}\\
0_{2mN}&=-\eta^*-(L\otimes I_{2m})\eta^*-(L\otimes I_{2m})w^*+\Theta, \label{cf12} \\
0_{2mN}&=(L\otimes I_{2m})\eta^*. \label{cf13}
\end{align}
\end{subequations}
Since the communication graph is undirected and connected, it follows from \eqref{cf13} that $\eta_1^*=\cdots=\eta_N^*$. According to the property of the undirected and connected graph $\mathcal{G}$, left multiplying \eqref{cf12} by $1_{2Nm}^T$ yields  $\eta_1^*=\cdots=\eta_N^*=(1/N)\sum_{j=1}^N\Theta_j$, which indicates that $\eta_{11}^*=\cdots=\eta_{N1}^*=(1/N)\sum_{j=1}^N\phi_j(x_j^*)$ and $\eta_{12}^*=\cdots=\eta_{N2}^*=(1/N)\sum_{j=1}^N\frac{\partial f_j}{\partial \sigma}$. By the definitions of $\eta_1$ and $\eta_2$, and \eqref{cf11}, we have that $\frac{\partial f_i(x_i^*,\sigma(x^*))}{\partial x_i}+(\nabla \phi_i(x_i^*)/N)\sum_{j=1}^N\frac{\partial f_j}{\partial \sigma}=0_{n_i}$ for all $i\in\mathcal{I}$,  which indicates that the optimal condition holds. Thus, $x^*$ is an optimal solution of problem \eqref{op}.

\textit{Sufficiency}: If $x^*$ is the optimal solution to problem \eqref{op}, it satisfies $\partial f_i(x_i^*,\sigma(x^*))/\partial x_i+(\nabla \phi_i(x_i)/N)(\sum_{j=1}^N\partial f_j(x_j^*,\sigma(x^*))/\partial \sigma)$$=$$0_{n_i}$ for all $i\in\mathcal{I}$, which indicates that $\eta_{11}^*=\cdots=\eta_{N1}^*=(1/N)\sum_{j=1}^N\phi_j(x_j^*)$ and $\eta_{12}^*=\cdots=\eta_{N2}^*=(1/N)\sum_{j=1}^N\frac{\partial f_j}{\partial \sigma}$. It follows from the definition of $\eta_i$ that $\eta_1^*=\cdots=\eta_N^*=(1/N)\sum_{j=1}^N\Theta_j$. Then, \eqref{cf11} and \eqref{cf13} hold. There exists $w^*$ such that \eqref{cf12} is satisfied. Thus, $(x^*,\eta^*,w^*)$ is an equilibrium of system \eqref{cf}.
$\hfill\blacksquare$

\begin{thm} \label{thm2}
The discrete-time communication among agents is triggered by rule \eqref{tc}, in which $\beta_{i1}>0$ and $0<\beta_{i2}<\lambda$ with $\lambda$ being the smallest positive eigenvalue of matrix $\big [\begin{smallmatrix}I+L & L \\ -L & 0\end{smallmatrix}\big ]$, for all $i\in\mathcal{I}$. Under Assumptions \ref{as2}-\ref{as3}, there exists a positive constant $\delta^*$ such that for every $0<\delta<\delta^*$, $\eta(t)$ converges exponentially to $(1/N)\sum_{i=1}^N \Theta_i 1_{N}$. The decision $x_i$ of agent $i$, $i\in\mathcal{I}$, which follows the designed algorithm \eqref{alg}-\eqref{dacd}, exponentially converges to the optimal decision $x_i^*$ associated with problem \eqref{op}. Moreover, the event-triggered communication scheme \eqref{tc} is free of the Zeno behavior.
\end{thm}

\textbf{Proof.} System \eqref{cf} can be considered as a singular perturbation system with parameter $\delta$. According to the stability analysis of the singular perturbation systems \cite{Khalil.2002}, the following analysis is divided into three steps.

1) Boundary-layer analysis: In the $\tau=t/\delta$ time scale, system \eqref{cf} is represented by
\begin{equation} \label{cft}
\begin{aligned}
\frac{dx}{d\tau}&=-\delta (F(x,\eta_1)+\eta_2\nabla \phi(x)),\\
\frac{d\eta}{d\tau}&=-\eta-(L\otimes I_{2m})\eta-(L\otimes I_{2m})w+\Theta\\
&\ \ \ -(L\otimes I_{2m})(E_1+E_2),\\
\frac{dw}{d\tau}&=(L\otimes I_{2m})\eta+(L\otimes I_{2m})E_1.
\end{aligned}
\end{equation}
Define $\tilde{\eta}=\eta-\bar{\eta}$ and $\tilde{w}=w-\bar{w}$ with steady states $\bar{\eta}$ and $\bar{w}$ of system \eqref{dacd}. Since $x$ in the foregoing equation is slowly varying, set $\delta=0$ to freeze $x$, and to reduce \eqref{cft} to the following system, called the boundary-layer system.
\begin{equation} \label{bsd1}
\begin{aligned}
\begin{bmatrix}
\frac{d \tilde{\eta}}{d\tau}\\ \frac{d \tilde{w}}{d\tau}
\end{bmatrix}
&=
\begin{bmatrix}
-I-(L\otimes I_{2m}) & -(L\otimes I_{2m}) \\
(L\otimes I_{2m}) & 0
\end{bmatrix}
\begin{bmatrix}
\tilde{\eta} \\ \tilde{w}
\end{bmatrix}
\\
&\ \ \ +
\begin{bmatrix}
-(L\otimes I_{2m})& -(L\otimes I_{2m}) \\ (L\otimes I_{2m}) & 0
\end{bmatrix}
\begin{bmatrix}
E_1 \\ E_2
\end{bmatrix}.
\end{aligned}
\end{equation}
Like the analysis procedure in \cite{Cai.2022,Ye.2017}, let $r$ be an $N$-dimensional column vector such that $r^TL=0$. Define $\Pi=[R,r]\otimes I_{2m}$ as an orthogonal matrix. Then, $\tilde{w}$ can be decomposed by $\tilde{w}=\Pi\begin{bmatrix} \tilde{w}' & \tilde{w}_o\end{bmatrix}^T$, in which $\tilde{w}_o$ is a vector in the consensus subspace and $\tilde{w}'$ is a vector in the orthogonal complement of the consensus subspace. Since $\tilde{w}_o$ does not interact with the other states, system \eqref{bsd1} can be reduced to
\begin{equation} \label{bsd2}
\begin{aligned}
\begin{bmatrix}
\frac{d \tilde{\eta}}{d\tau}\\ \frac{d \tilde{w}'}{d\tau}
\end{bmatrix}
&=
\underbrace{\begin{bmatrix}
-I-(L\otimes I_{2m}) & -LR\otimes I_{2m} \\
R^TL\otimes I_{2m} & 0
\end{bmatrix}}_{P}
\begin{bmatrix}
\tilde{\eta} \\ \tilde{w}'
\end{bmatrix}
\\
&\ \ \ +
\underbrace{\begin{bmatrix}
-(L\otimes I_{2m}) & -(L\otimes I_{2m}) \\ (R^TL\otimes I_{2m}) & 0
\end{bmatrix}}_{Q}
\begin{bmatrix}
E_1 \\ E_2'
\end{bmatrix}.
\end{aligned}
\end{equation}
It is easily known that system \eqref{bsd1} converges to the origin if system \eqref{bsd2} converges to the origin.
Since the communication graph is undirected and connected, $I+L>0$ and $R^TL\otimes I_{2m}$ is of full rank. It follows from Lemma 2.2 in \cite{Menon.2014} that matrix $P$ is Hurwitz. It is easily derived that $\|[\tilde{\eta}, \tilde{w}']\|\leq e^{-\lambda \tau}\|[\tilde{\eta}(0), \tilde{w}'(0)]\|+\int_0^\tau e^{-\lambda(\tau-s)}\|Q\| \|\col(E_1,E_2')\|ds \leq e^{-\lambda \tau}\|[\tilde{\eta}(0), \tilde{w}'(0)]\|+\int_0^\tau e^{-\lambda(\tau-s)}\|Q\| \|\col(E_1,E_2)\|ds \leq e^{-\lambda \tau}\|[\tilde{\eta}(0), \tilde{w}'(0)]\|+\beta_1\sqrt{N}\|Q\|(e^{-\lambda \tau}-e^{-\beta_2 \tau})/(\beta_2-\lambda)$,
where the third inequality comes from the event-triggering condition \eqref{tc}, $\beta_1$$=$$\max\{\beta_{i1},i=1,\ldots,N\}$, $\beta_2$$=$$\min\{\beta_{i2},i=1,\ldots,N\}$, and $\lambda$ is the absolute value of the largest eigenvalue of matrix $P$. Therefore, system \eqref{bsd2} exponentially converges to the origin. It indicates that the dynamic average consensus estimator \eqref{dacd} with discrete-time communication can exponentially converge to $(1/N)\sum_{j=1}^N\Theta_j$.

2) Quasi-steady state analysis: Set $\delta=0$ such that $\eta$ and $w$ are at their quasi-steady states $\bar{\eta}$ and $\bar{w}$, respectively, in which $\bar{\eta}_{i1}=(1/N)\sum_{j=1}^N\phi_j(x_j)$ and $\bar{\eta}_{i2}=(1/N)\sum_{j=1}^N\partial f_j/\partial \sigma$, $\forall i\in\mathcal{I}$. Then, system \eqref{cf} can be reduced to
\begin{equation}
\begin{aligned} \label{rs}
\dot{x}=-F_1(x,\sigma)-F_2\nabla \phi(x)
=-\nabla f(x),
\end{aligned}
\end{equation}
where $F_1(x,\sigma)$ $=$ $\col(\partial f_1(x_1,\sigma)/\partial x_1,$ $\ldots,$ $\partial f_N(x_N,\sigma)/\partial x_N)$ and $F_2$ $=$ $(1/N)\sum_{i=1}^N\partial f_i/\partial \sigma 1_{N}$, thus the second equality is easily obtained. It follows from Lemma \ref{thm1} that the reduced system \eqref{rs} exponentially converges to the optimal decision $x^*$ at a convergence rate no less than $e^{-\kappa/(1+2l)t}$.

3) Convergence analysis: According to Theorem 11.4 in \cite{Khalil.2002}, it is concluded that there exists a positive constant $\delta^*>0$ such that for all $\delta\in(0,\delta^*)$, all agents' decisions converge exponentially to the optimal solution $x^*$ of problem \eqref{op}, i.e., $\lim_{t\rightarrow \infty}\|x-x^*\|=0$. In other words, each agent's decision exponentially converges to its optimal decision $x_i^*$ in the global sense.

%The analysis of the quasi-steady state is similar to that in the proof of Theorem \ref{thm1} and is omitted due to the limited space. By Theorem 11.4 in \cite{Khalil.2002}, it is concluded that there exists a positive constant $\delta^*$ such that for all $\delta\in(0,\delta^*)$, all agents' decisions can converge exponentially to the optimal solution $x^*$ of problem \eqref{op}, if the time instants when agents' information is broadcasted are determined by the triggering rule \eqref{tc}.

Next, the Zeno behavior in the triggering rule \eqref{tc} is analyzed. For agent $i$, the dynamics of $e_i$ is given by $\dot{e}_i=\col(\dot{e}_{i1},\dot{e}_{i2})$ with $\dot{e}_{i1}=\dot{\eta}_{i}$ and $\dot{e}_{i2}=\dot{w}_i$. It is derived that $d\|e_i\|/dt\leq \|de_i/dt\|\leq \|\col(\dot{\eta}_i,\dot{w}_i)\|\leq \|\col(\dot{\eta},\dot{w})\|\leq \|\col(\dot{\tilde{\eta}},\dot{\tilde{w}})\|\leq \|\col(\dot{\tilde{\eta}},\dot{\tilde{w}}')\|$.
It follows from the analysis of boundary-layer system \eqref{bsd1} that $\|\col(\dot{\tilde{\eta}},\dot{\tilde{w}}')\|\leq \|P\| \|\col(\tilde{\eta},\tilde{w}')\|+\|Q\| \|\col(E_1,E_2')\|\leq M_1 e^{-\lambda t}+M_2e^{-\beta_2 t}$ with $M_1=\|P\| \|\col(\tilde{\eta}(0),\tilde{w}'(0))\|+\sqrt{N}\beta_1\|P\|\|Q\|/(\beta_2-\lambda)$ and $M_2=\|Q\|\beta_1\sqrt{N}$.
Assume that $t^*$ is the latest triggering instant. The upper bound of $\|\col(\dot{\tilde{\eta}},\dot{\tilde{w}}')\|$, depending on $t^*$, is $M_1e^{-\lambda t^*}+M_2e^{-\beta_2 t^*}$. Thus, it yields that $\|e_i(t)\|\leq (M_1e^{-\lambda t^*}+M_2e^{-\beta_2 t^*})(t-t^*)$.
The next triggering instant is the time when the inequality $\|e_i\|\geq\beta_{i1}e^{-\beta_{i2}t}$ holds.
Let $T$ denote a lower bound of inter-event intervals, which is the solution of $(M_1e^{(\beta_2-\lambda)t^*}+M_2)T=\beta_{i1}e^{-\beta_{i2}T}$. For $\beta_2<\lambda$, (i.e.,  $\beta_{i2}<\lambda$ for all $i\in\mathcal{I}$), and $t^*\geq0$, $T(t^*)$ is greater or equal to $T$ given by $(M_1+M_2)T=\beta_{i1}e^{-\beta_{i2}T}$, which is strictly positive. It implies that there exists a positive lower bound on inter-event intervals. Therefore, the Zeno behavior is excluded.
$\hfill\blacksquare$

\begin{rmk}
Singular perturbation systems have a multitime-scale characteristic. To be specific, in system \eqref{cf}, state $x$ is slowly varying and state $(\eta,w)$ is varying relatively fast. In the application of high-gain feedback systems, parameter $\delta$ is the reciprocal of a high-gain parameter by simple transformation of equation \eqref{cf}. The singular perturbation method is a typical approach to analyzing and designing high-gain feedback systems.
\end{rmk}

\section{Case studies} \label{sec4}
\begin{figure}[!t]
  \centering
  \includegraphics[width=4cm]{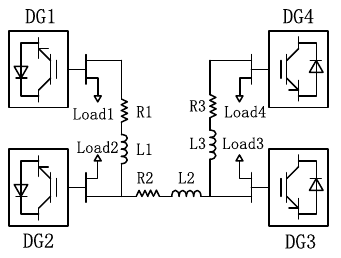}
  \caption{The structure of the four DERs.}
  \label{fig1}
\end{figure}

\begin{figure}
\begin{center}
\subfigure[]
{\includegraphics[width=4cm,height=3.6cm]{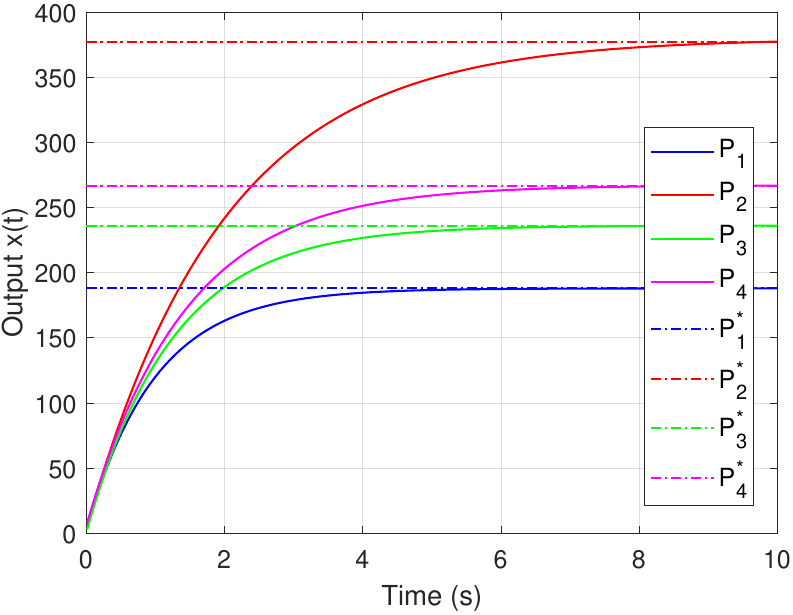}}
\subfigure[]
{\includegraphics[width=4cm,height=3.6cm]{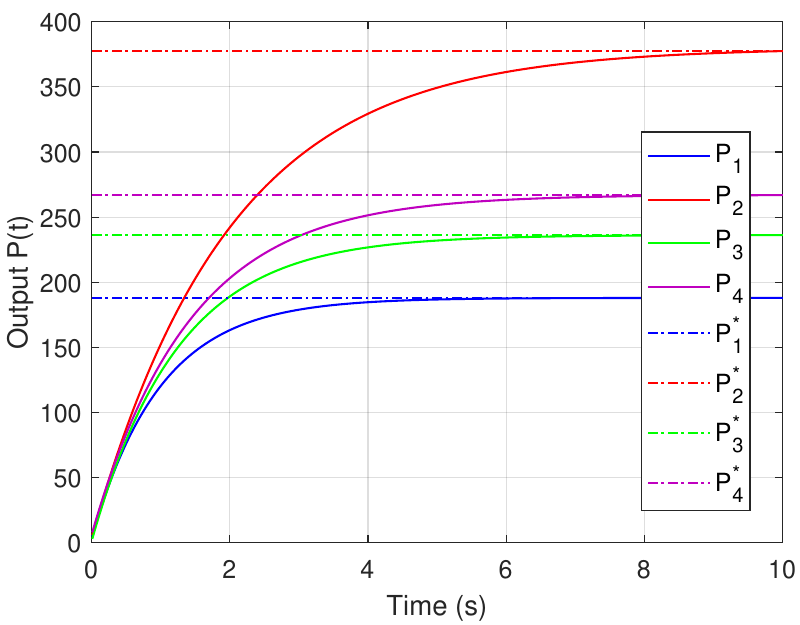}}
\DeclareGraphicsExtensions.
\caption{Evolution of DERs' decisions: (a) continuous-time communication (b) event-triggered communication.}
\label{fig3}
\end{center}
\end{figure}

\begin{figure}
\begin{center}
\subfigure[]
{\includegraphics[width=4cm,height=3.6cm]{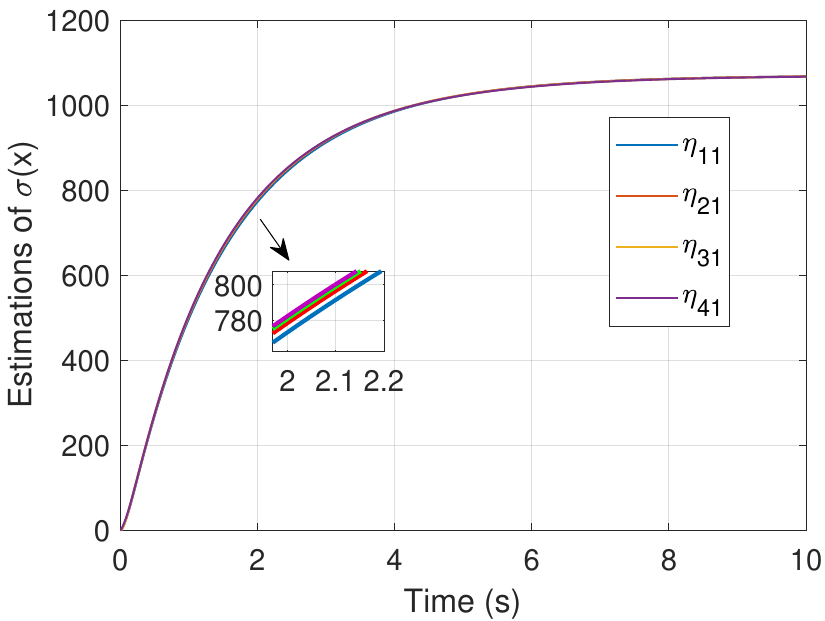}}
\subfigure[]
{\includegraphics[width=4cm,height=3.6cm]{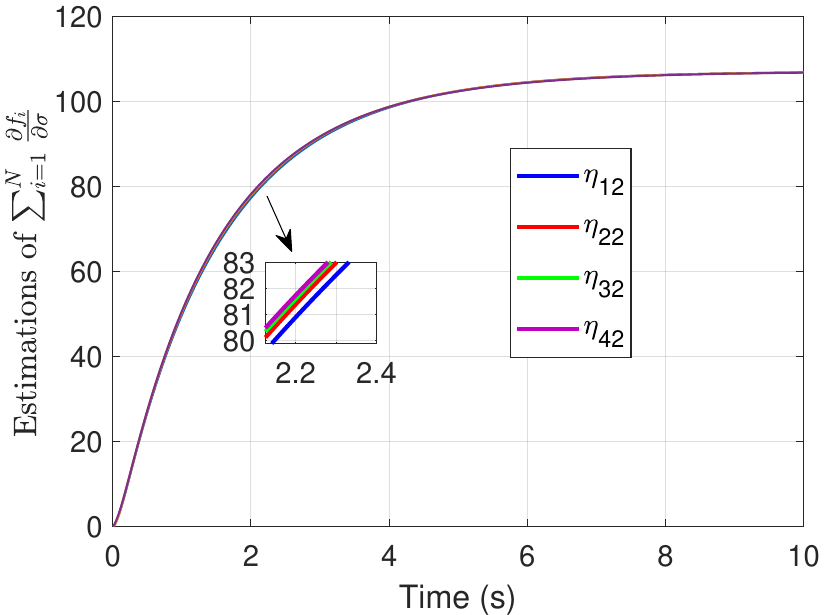}}
\DeclareGraphicsExtensions.
\caption{The output of estimator \eqref{dacd}: (a) estimations of aggregator $\sigma$ (b) estimations of $(1/N)\sum_{i=1}^N \partial f_i/\partial \sigma$.}
\label{fig6}
\end{center}
\end{figure}

\begin{figure}[!t]
  \centering
  \includegraphics[width=6cm]{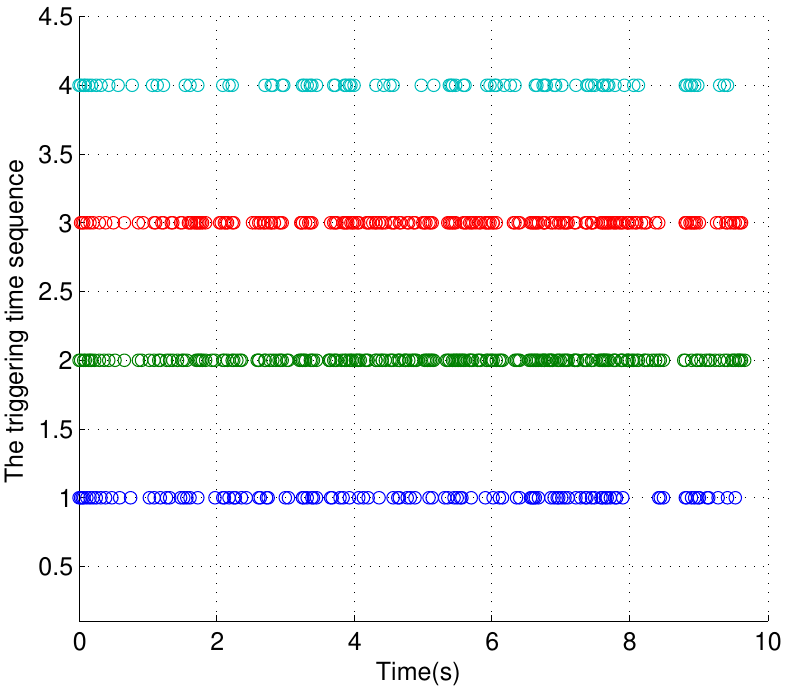}
  \caption{Event-triggering time sequences.}
  \label{fig4}
\end{figure}

\begin{figure}[!t]
  \centering
  \includegraphics[width=6cm]{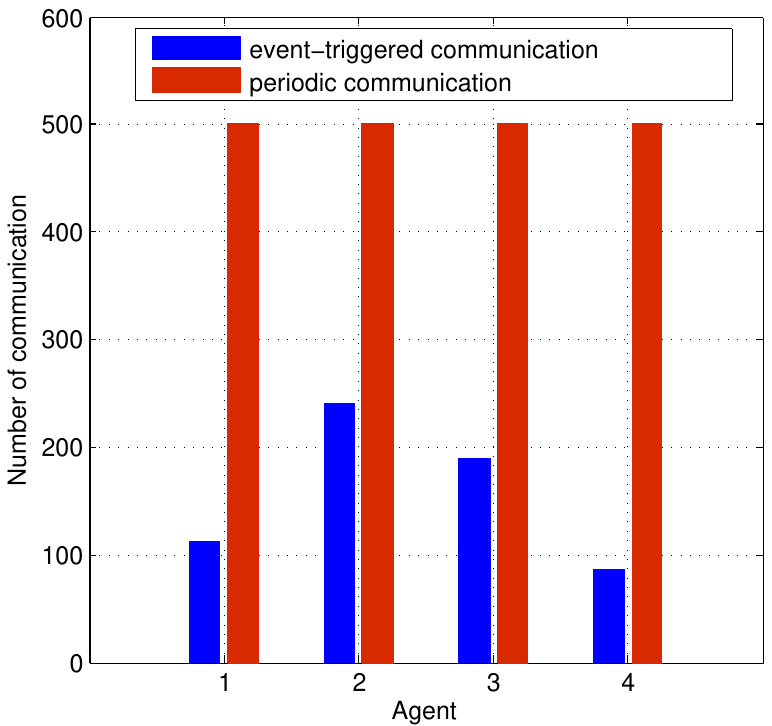}
  \caption{Comparison between communication times of event-triggering rule \eqref{tc} and of period sampling $T=0.02s$.}
  \label{fig5}
\end{figure}

In this section, taken the price-based energy management of DERs as an example, the four DERs, forming a microgrid (referring to \cite{Chen.2022}), are shown in Fig.~\ref{fig1}.
In the simulation, the initial conditions are $x(0)=[5, 6, 3, 8]^T$, $a=[1.0 0.5 0.8 0.7]$, $b=[12, 10, 11, 11]$, $d=[5, 8, 6, 9]$, and $\delta=0.1$. The optimal solution is $P^*=[188.0, 377.5, 236.2, 266.9]^T$ by simple calculation. In event-triggered communication scheme \eqref{tc}, $[\beta_{11},\ldots,\beta_{14}]^T=[10, 8, 8, 10]^T$ and $[\beta_{21},\ldots,\beta_{24}]^T=[0.01, 0.1, 0.15, 0.05]^T$. Fig.~\ref{fig3} depicts all DERs' decisions obtained by the designed algorithm \eqref{alg}-\eqref{dacd} in the continuous-time communication and the event-triggered broadcasting scheme \eqref{tc}, respectively. Fig.~\ref{fig3} indicates that all DERs' decisions can reach their own optimal decision $P_i^* (i\in\{1,\ldots,4\})$. Moreover, it is seen from Fig.~\ref{fig3} that the proposed event-triggered communication scheme has little influence on the convergence of algorithm \eqref{alg}. The outputs of the dynamic average estimator \eqref{dacd} is shown in Fig.~\ref{fig6}, which indicates that the estimator \eqref{dacd} can quickly estimate the global information.  Fig.~\ref{fig4} shows the triggering instants of the four DERs. It can be seen that the Zeno behavior is excluded. It can be seen from Fig.~\ref{fig5} that the total number of communications of event-triggering rule \eqref{tc} are less than that of the period sampling with $T=0.02s$. It implies that the designed event-triggered broadcasting scheme can reduce communication loads effectively.

\begin{figure}[!t]
  \centering
  \includegraphics[width=6cm]{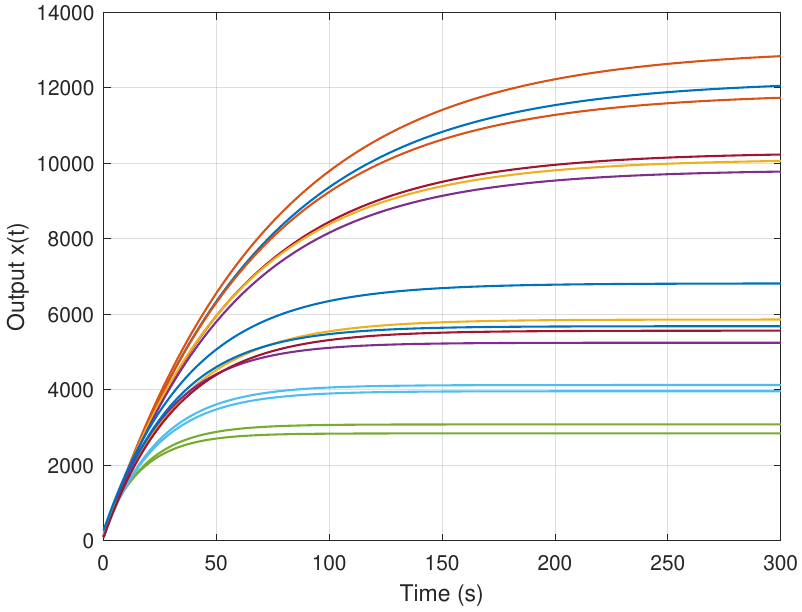}
  \caption{Outputs of the fifteen generators.}
  \label{fig2}
\end{figure}

Next, to show the effectiveness of the proposed distributed algorithm for large-scale systems, consider an IEEE 118-bus system with fifteen generators \cite{Yi.2016} as an example. The communication network is a randomly generated undirected and connected graph. Parameters $a_i$, $b_i$, and $d_i$ in cost function $f_i$ belong to the intervals $a_i\in[0.0024,0.0779]$, $b_i\in[8,35]$, and $d_i\in[7, 60]$, respectively. Set $\delta=0.1$, $\beta_{i1}=6$, and $\beta_{i2}=0.15$ in \eqref{dacd}-\eqref{tc} for $i\in\{1,2,\ldots,15\}$. The outputs of the fifteen generators are shown in Fig.~\ref{fig2}, from which the outputs of the fifteen generators reach to their optimal decisions.

\section{Conclusion} \label{sec5}
In this paper, a distributed continuous-time optimization algorithm with event-triggered communication has been proposed for solving a distributed aggregative optimization problem with constrained communication, in which local optimization functions depend on agents' own decisions and the aggregation of all agents' decisions. Rigorous analysis has been given to ensure the exponential convergence of the designed algorithm.  A numerical example has been provided to illustrate the obtained results. The interplay between the event-triggered communication times and the convergence rate of the algorithm will be studied as an interesting issue in the future.

\appendix

\bibliographystyle{plain}        % Include this if you use bibtex
\bibliography{reference}

\begin{thebibliography}{10}

\bibitem{Boyd.2004}
S.~Boyd and L.~Vandenberghe.
\newblock {\em Convex optimization}.
\newblock Cambridge University Press, USA, 2004.

\bibitem{Cai.2021}
X.~Cai, F.~Xiao, and B.~Wei.
\newblock A distributed event-triggered generalized {N}ash equilibrium seeking
  algorithm.
\newblock In {\em 2021 40th Chinese Control Conference (CCC)}, pages
  5252--5257, 2021.

\bibitem{Cai.2022}
X.~Cai, F.~Xiao, and B.~Wei.
\newblock Nash equilibrium seeking for general linear systems with disturbance
  rejection.
\newblock {\em IEEE Transactions on Cybernetics}, 53(8):5240--5249, 2023.

\bibitem{Cai.2021a}
X.~Cai, F.~Xiao, B.~Wei, and F.~Fang.
\newblock Distributed continuous-time strategy-updating rules for
  noncooperative games with discrete-time communication.
\newblock {\em IEEE Transactions on Systems, Man, and Cybernetics: Systems},
  pages 1--10, 2023.

\bibitem{Carnevale.2022a}
G.~Carnevale, A.~Camisa, and G.~Notarstefano.
\newblock Distributed online aggregative optimization for dynamic multi-robot
  coordination.
\newblock {\em IEEE Transactions on Automatic Control}, pages 1--8, 2022.

\bibitem{Carnevale.2022}
G.~Carnevale, F.~Fabiani, F.~Fele, K.~Margellos, and G.~Notarstefano.
\newblock Tracking-based distributed equilibrium seeking for aggregative games.
\newblock {\em arXiv:2210.14547v1}, pages 1--41, 2022.

\bibitem{Chen.2022}
J.~Chen, D.~Yue, C.~Dou, S.~Weng, X.~Xie, Y.~Li, and G.~Hancke.
\newblock Static and dynamic event-triggered mechanisms for distributed
  secondary control of inverters in low-voltage islanded microgrids.
\newblock {\em IEEE Transactions on Cybernetics}, 52(7):6925--6938, 2022.

\bibitem{Chen.2016}
W.~Chen and W.~Ren.
\newblock Event-triggered zero-gradient-sum distributed consensus optimization
  over directed networks.
\newblock {\em Automatica}, 65:90--97, 2016.

\bibitem{Dai.2020}
H.~Dai, X.~Fang, and W.~Chen.
\newblock Distributed event-triggered algorithms for a class of convex
  optimization problems over directed networks.
\newblock {\em Automatica}, 122:109256, 2020.

\bibitem{Gharesifard.2016}
B.~Gharesifard, T.~Ba\c{s}ar, and A.~Dom\'{\i}nguez-Garc\'{\i}a.
\newblock Price-based coordinated aggregation of networked distributed energy
  resources.
\newblock {\em IEEE Transactions on Automatic Control}, 61(10):2936--2946,
  2016.

\bibitem{Grammatico.2017}
S.~Grammatico.
\newblock Dynamic control of agents playing aggregative games with coupling
  constraints.
\newblock {\em IEEE Transactions on Automatic Control}, 62(9):4537--4548, 2017.

\bibitem{Khalil.2002}
H.K. Khalil.
\newblock {\em Nonlinear Systems}.
\newblock Prentice-Hall, USA, 3rd edition, 2002.

\bibitem{Koshal.2016}
J.~Koshal, A.~Nedi\'{c}, and U.~Shanbhag.
\newblock Distributed algorithms for aggregative games on graphs.
\newblock {\em Operations Research}, 64(3):680--704, 2016.

\bibitem{Li.2021}
X.~Li, Xie L, and Y.~Hong.
\newblock Distributed aggregative optimization over multi-agent networks.
\newblock {\em IEEE Transactions on Automatic Control}, 67(6):3165--3171, 2022.

\bibitem{Liang.2019}
S.~Liang, L.~Y. Wang, and G.~Yin.
\newblock Exponential convergence of distributed primal-dual convex
  optimization algorithm without strong convexity.
\newblock {\em Automatica}, 105:298--306, 2019.

\bibitem{Liang.2017}
S.~Liang, P.~Yi, and Y.~Hong.
\newblock Distributed {N}ash equilibrium seeking for aggregative games with
  coupled constraints.
\newblock {\em Automatica}, 85:179--185, 2017.

\bibitem{Liang.2022}
S.~Liang, P.~Yi, Y.~Hong, and K.~Peng.
\newblock Exponentially convergent distributed {N}ash equilibrium seeking for
  constrained aggregative games.
\newblock {\em Autonomous Intelligent Systems}, 2:6, 2022.

\bibitem{Liu.2020}
C.~Liu, H.~Li, Y.~Shi, and D.~Xu.
\newblock Distributed event-triggered gradient method for constrained convex
  minimization.
\newblock {\em IEEE Transactions on Automatic Control}, 65(2):778--785, 2020.

\bibitem{Menon.2014}
A.~Menon and J.~Baras.
\newblock Collaborative extremum seeking for welfare optimization.
\newblock In {\em 53th IEEE Conference on Desicion and Control}, pages
  345--351, 2014.

\bibitem{Romano.2020}
A.~R. Romano and L.~Pavel.
\newblock Dynamic {NE} seeking for multi-integrator networked agents with
  disturbance rejection.
\newblock {\em IEEE Transactions on Control of Network Systems}, 7(1):129--139,
  2020.

\bibitem{Seyboth.2013}
G.~S. Seyboth, D.~V. Dimarogonas, and K.~H. Johansson.
\newblock Event-based broadcasting for multi-agent average consensus.
\newblock {\em Automatica}, 49:245--252, 2013.

\bibitem{Shakarami.2022}
M.~Shakarami, C.~De Persis, and N.~Monshizadeh.
\newblock Distributed dynamics for aggregative games: Robustness and privacy
  guarantees.
\newblock {\em International Journal of Robust and Nonlinear Control},
  32:5048--5069, 2022.

\bibitem{Wang.2019a}
A.~Wang, B.~Mu, and Y.~Shi.
\newblock Event-triggered consensus control for multiagent systems with
  time-varying communication and event-detecting delays.
\newblock {\em IEEE Transactions on Control Systems Technology},
  27(2):507--515, 2019.

\bibitem{Wang.2022}
T.~Wang and P.~Yi.
\newblock Distributed projection-free algorithm for constrained aggregative
  optimization.
\newblock {\em arXiv:2207.11885v1}, pages 1--21, 2022.

\bibitem{Ye.2017}
M.~Ye and G.~Hu.
\newblock Game design and analysis for price-based demand response: An
  aggregate game approach.
\newblock {\em IEEE Transactions on Cybernetics}, 47(3):720--730, 2017.

\bibitem{Ye.2021}
M.~Ye, G.~Hu, L.~Xie, and S.~Xu.
\newblock Differentially private distributed {N}ash equilibrium seeking for
  aggregative games.
\newblock {\em IEEE Transactions on Automatic Control}, 67(5):2451--2458, 2022.

\bibitem{Yi.2016}
P.~Yi, Y.~Hong, and F.~Liu.
\newblock Initialization-free distributed algorithms for optimal resource
  allocation with feasibility constraints and application to economic dispatch
  of power systems.
\newblock {\em Automatica}, 74:259--269, 2016.

\bibitem{Yi.2020}
X.~Yi, S.~Zhang, T.~Yang, T.~Chai, and K.~H. Johansson.
\newblock Linear convergence for distributed optimization without strong
  convexity.
\newblock In {\em 59th IEEE Conference on Decision and Control}, pages
  3643--3648, 2020.

\bibitem{Yu.2021}
H.~Yu and T.~Chen.
\newblock A new {Z}eno-free event-triggered scheme for robust distributed
  optimal coordination.
\newblock {\em Automatica}, 129:109639, 2021.

\end{thebibliography}
%\section{A summary of Latin grammar}    % Each appendix must have a short title.
%\section{Some Latin vocabulary}         % Sections and subsections are supported
                                        % in the appendices.
\end{document}